\documentclass[a4paper,12pt]{article}

\usepackage{amssymb}

\newtheorem{lemma}{Lemma}
\newtheorem{proposition}[lemma]{Proposition}
\newtheorem{theorem}[lemma]{Theorem}
\newtheorem{corollary}[lemma]{Corollary}

\newtheorem{conjecture}[lemma]{Conjecture}

\newcommand{\alp}{\alpha}
\newcommand{\bet}{\beta}
\newcommand{\gam}{\gamma}
\newcommand{\kap}{\kappa}
\newcommand{\lam}{\lambda}

\newcommand{\del}{\delta}
\newcommand{\eps}{\varepsilon}
\newcommand{\tet}{\vartheta}

\newcommand{\Ome}{\Omega}

\newcommand{\BB}{{\cal B}}
\renewcommand{\Re}{{\rm Re}}
\renewcommand{\Im}{{\rm Im}}

\newcommand{\R}{\mathbb{R}}
\newcommand{\C}{\mathbb{C}}
\newcommand{\Lt}{{\rm L}^2}
\newcommand{\Ci}{{\rm C}^{\infty}}
\newcommand{\Ic}{{\rm C}^{\infty}_{\rm c}}
\newcommand{\Sc}{{\cal S}(\R)}
\newcommand{\So}{{\rm W}^{1,2}}

\newcommand{\Schrodinger}{Schr\"o\-din\-ger }
\newcommand{\bel}{ \in }
\newcommand{\e}{{\rm e}}
\newcommand{\inte}{\mathrm{Int}}

\newcommand{\norm}{\Vert}
\newcommand{\spec}[1]{{\rm Spec} \left( #1 \right) \, }
\newcommand{\pspe}[2]{{\rm Spec} _{#1} \left( #2 \right) \, }

\newcommand{\nume}[1]{{\rm Num} \left( #1 \right) \, }

\newcommand{\srad}[1]{{\rm rad} \left( #1 \right) \, }

\newcommand{\rang}[1]{{\rm Ran} \left( #1 \right) \, }

\newcommand{\gap}{\vspace{.2in}}
\newcommand{\Proof}{\underline{Proof} {\hskip 0.1in}}
\newcommand{\fin}{\hspace{.1in} $\, \blacksquare$}
\newcommand{\find}{\hspace{.1in} \blacksquare}

\title{NON-SELF-ADJOINT HARMONIC OSCILLATOR, \\
     COMPACT SEMIGROUPS AND PSEUDOSPECTRA}
\author{Lyonell S. Boulton}
\date{29 September 1999}

\begin{document}

\maketitle

\begin{abstract}
\setlength{\baselineskip}{15pt} We provide new information
concerning the pseudospectra of the complex harmonic oscillator.
Our analysis illustrates two different techniques for getting
resolvent norm estimates. The first uses the JWKB method and
extends for this particular potential some results obtained
recently by E.B. Davies. The second relies on the fact that the
bounded holomorphic semigroup generated by the complex harmonic
oscillator is of Hilbert-Schmidt type in a maximal angular region.
In order to show this last property, we deduce a non-self-adjoint
version of the classical Mehler's formula.
\end{abstract}

\setlength{\parindent}{0cm}


\vskip.1in


\setlength{\parindent}{.6cm}

\pagebreak

\section{Introduction}

The harmonic oscillator is known to be ubiquitous in many branches
of quantum theory. Unlike most quantum mechanical problems, it can
be solved explicitly in terms of elementary functions (see section
2), so it provides an excellent illustration for the general
principles of the quantum-theoretical formalism. It also plays a
key role in quantum electrodynamics and quantum field theory, and
has been a focus of great attention among the mathematical-physics
community for many years.

Although it might seem very difficult to say something new about
the harmonic oscillator, in a couple of recent papers
(\cite{D2,D3}) E. B. Davies provides a framework which opens the
possibility of obtaining new information concerning the spectral
theory of this and other \Schrodinger operators. In the present
paper we extend some results given in \cite{D2} and \cite{D3}, and
introduce some new ideas that will be useful in studying both
spectral and stability properties of this operator when it has a
complex coupling constant.

We define the complex harmonic oscillator to be the operator
\[
   H_cf(x):=-\frac{d^2}{dx^2}f(x)+cx^2f(x)
\]
acting on $\Lt(\R)$ with Dirichlet boundary conditions, where
$c\bel \C$ is such that $\Re(c)>0$ and $\Im(c)>0$. Our main
interest is to investigate the resolvent norm of $H_c$ inside its
numerical range. To this end, if we let the $\eps$-pseudospectra
of $H_c$ be
\[
   \pspe{\eps}{H_c} := \spec{H_c} \cup \{ z \bel \C \,:\, \norm (H_c-z)^{-1}
   \norm \geq \frac{1}{\eps} \}
\]
for all $\eps>0$, our interest is in obtaining information
concerning the shape of these sets for small $\eps$.

We adopt the notation and basic ideas about pseudospectra recently
developed by E. B. Davies, L. N. Trefethen and others. For general
results and more examples in the present spirit, we refer to
\cite{D1,D2,T1} and to the bibliography there.

The results we obtain here extend in two ways those given by E. B.
Davies in \cite{D2,D3} where he shows that for any positive
constant $b$ independent of $\eta>0$,
\[
   \norm (H_c-\eta[b+c])^{-1} \norm \to \infty
\]
as $\eta \to \infty$. On the one hand in section 3 (theorems
\ref{t18}) we demonstrate that for all $b>0$ and $1/3 < p < 3$
fixed,
\[
   \norm \left( H_c-[b\eta+c\eta^p] \right)^{-1} \norm \to \infty
\]
as $\eta\to \infty$. On the other hand in section 5 we prove (as a
consequence of theorem \ref{t15}) that for all $b>0$ there exist a
constant $M_b>0$ such that
\[
   \lim_{\eta \to \infty} \norm (H_c-[\eta+ib])^{-1} \norm \leq M_b
\]
and
\[
   \lim_{\eta \to \infty} \norm (H_c-[c \eta-i b])^{-1} \norm  \leq
   M_b.
\]

The study of pseudospectra provides important information about
the stability of operator $H_c$. By virtue of the identity
\[
   \pspe{\eps}{H_c} = \bigcup \{ \spec{H_c+A} \,:\, \norm A \norm \leq
   \eps\}
\]
(for a proof see \cite{RS1}), knowing the size of these sets for
$\eps$ close to zero, allow us to obtain very precise information
about the stability under small perturbation of the eigenvalues of
$H_c$, \cite{D1}. It will turn out (sections 3 and 5) that high
energy eigenvalues are far more unstable than the first excited
states (see also \cite{AD1,D2}).

The method we use in section 3 is analogous to the techniques in
\cite{D2,D3} and involves the construction of a continuous family
of approximate eigenstates for $H_c$ by means of JWKB analysis
(see \cite{Hel} for an introduction to this topic and some
applications). We should mention that such a procedure is examined
for more general potentials and from the numerical point of view
in \cite{AD2}.

In section 4 we deduce a non-self-adjoint version of the classical
Mehler's formula, {\it i.e.} we construct an explicit formula for
the heat kernel of $-H_c$ (theorem \nolinebreak \ref{t3}). This
formula allows us to show that the bounded holomorphic semigroup
generated by $-H_c$ is compact in a maximal angular domain. The
results of section 4 provides a new approach to get resolvent norm
estimates for $H_c$. Based on these estimates, in section 5 we
obtain a set that encloses $\pspe{\eps}{H_c}$ for small enough
$\eps$ (corollary \ref{t16} and theorem \ref{t15}). This confirms
the numerical evidence given in \cite{D2} about the shape of such
sets.

\section{Definitions and Notation}

We will suppose in this paper that the parameter $c\bel \C$
satisfies $\Re(c)>0$ and unless explicitly stated we do not impose
conditions under $\Im(c)$.

We assume that $H_c$ acts in $\Lt(\R)$ with Dirichlet boundary
conditions as follows. Take the closed m-sectorial quadratic form
\[
   Q_c(f,g):=\int _\R f'(x)\overline{g'(x)} dx
   + c\int _\R x^2 f(x)\overline{g(x)} dx
\]
for all $f,g \bel \So(\R) \cap \{ f\bel \Lt(\R)  : \int x^2
|f(x)|^2 dx < \infty \}$. Then $H_c$ is defined to be the
m-sectorial operator associated to $Q_c$ via the Friedrichs
representation theorem (see \cite{kat}).

The subspaces $\Ic(\R)$ and $\Sc$ (the Schwartz space) are form
cores for $H_c$. Since
\[
    \Ic(\R) \subset \Sc \subset \So(\R) \cap \{ f\bel \Lt(\R)  : \int x^2
|f(x)|^2 dx < \infty \},
\]
it is enough to check the desired property for $\Ic(\R)$. By
\cite[theorem 1.21, p.317]{kat} a subspace is a core for an
m-sectorial form, if and only if it is a core for its real part.
Notice that the real part of $Q_c$ is the non-negative quadratic
form $\Re(Q_c)=Q_{\Re(c)}$, thus using for example \cite[theorem
1.13]{cyc} we can deduce that $\Ic(\R)$ is a core for $\Re(Q_c)$
and therefore that it is also a core for $Q_c$ as needed.

If $c$ is positive, $H_c$ is self-adjoint, but if $\Im(c)\neq 0$,
$H_c$ is not even a normal operator. Moreover (see \cite{D2})
there does not exist an invertible operator $U$ such that
$UH_cU^{-1}$ is normal.

If we put $\lam_n:=c^{1/2}(2n+1)$, the spectrum
of $H_c$ is
\[
   \spec{H_c}=\{\lam_n\,:\, n=0,\,1,\, \ldots \}.
\]
It consists entirely of eigenvalues of multiplicity one and it is
routine to check that if $H_n$ is the $n^{th}$-Hermite polynomial,
then
\[
   \Psi_n(x):=c^{1/8}H_n(c^{1/4}x)\e^{-(c^{1/2}x^2) / 2}
\]
are the eigenfunctions of $H_c$, so that
\[
   H_c \Psi_n = \lam_n \Psi_n \hspace{.3in} (n=0,\,1,\,\ldots).
\]

For $0\leq \alp,\,\bet \leq \pi/2$, we shall denote the angular
sector
\[
   S(-\alp,\bet):=\{z\bel \C \,:\, -\alp < \arg(z) < \bet
   \}.
\]
We will also put $S(\alp):=S(-\alp,\alp)$. It is clear from the
definition that  the numerical range of the operator $H_c$ is
contained in the angular sector $S(0,\arg(c))$.

\begin{proposition}
\label{t20}
For all $\Re(c)>0$,
\[
   \nume{H_c}= \{t_1+ct_2 \bel \C \, :\, t_1,\,t_2 \geq 0,\, t_1 t_2 \geq
   1/4 \}.
\]
If $\Im(c)\not= 0$, then
\[
   \spec{H_c}\subset \inte (\nume{H_c}).
\]
\end{proposition}

\Proof Since $\Sc$ is a core for $Q_c$,
\[
   \overline{\nume{H_c}}=\overline{ \left\{ Q_c(f) \, :
   \norm f \norm =1,\, f\bel \Sc \right\} }.
\]
By the Heisenberg inequality, we know that for any $f\bel \Sc$
\[
\left(\int _\R |f'(x)|^2 dx\right)
   \left(\int _\R x^2 |f(x)|^2 dx\right)
   \geq \frac{1}{4} \norm f \norm ^4 .
\]
Therefore clearly
\[
   \nume{H_c} \subset W:=\left\{ t_1+ct_2 \bel \C \, :
   \, t_1,t_2 > 0,\, t_1 t_2 \geq 1/4
   \right\} .
\]

Let us check the reverse inclusion. For this we need to find
test functions $f \bel \Sc$ such that
\[
   Q_c(f) \bel \partial (W)=\left\{ t+\frac{c}{4t} \bel \C \, :
   \, t > 0 \right\}.
\]
For $t>0$, let
\[
   f_{t} (x):= \e ^{-t x^2} \bel \Sc.
\]
Using elementary properties of the Gamma function we can
calculate
\begin{eqnarray*}
   \int _\R |f_{t} '(x)|^2 dx& = & 2t ^2 \pi ^{1/2} (2t)^{-3/2}, \\
   \int _\R x^2 |f_{t}(x)|^2 dx& = &  \frac{\pi ^{1/2}}{2} (2t)^{-3/2}, \\
   \norm f_{t} \norm ^2 & = & \pi^{1/2} (2t)^{-1/2}.
\end{eqnarray*}
Combining this three equalities, we obtain for all $t>0$
\[
\left(\int _\R |f'_{t} (x)|^2 dx\right)
   \left(\int _\R x^2 |f_{t}(x)|^2 dx\right)
   = \frac{1}{4} \norm f_{t} \norm ^4.
\]
Therefore
\[
   Q_c\left(\frac{f_{t}}{\norm f_{t} \norm}\right) = t+\frac{c}{4t}
   \hspace{.3in}(t >0)
\]
and so every point in $\partial (W)$ is in $\nume{H_c}$.
Since both sets are convex we have {\it a fortiori}
\[
   W\subset \nume{H_c}.
\]

For the second part, suppose that $\Im(c)\not=0$. Since
\[
   \spec{H_c} \subset \nume{H_c},
\]
it is enough to show that
\[
   \partial (W) \cap \spec{H_c}=\emptyset .
\]
The eigenvalues of $H_c$ lie on the line $c^{1/2}r,\, r>0$. The
boundary $\partial (W)$ only intersect this line at the point
\[
   z_0:=\frac{|c|^{1/2}}{2} \left(1+\frac{c}{|c|}\right).
\]
Clearly $|z_0|<|c^{1/2}|$, therefore the desired property follows
from the fact that $z_0$ is never an eigenvalue of $H_c$. \fin

 \gap

From the basic theory (see for instance \cite{D1,T1}), we know
that $\pspe{\eps}{H_c}$ contains the $\eps$-neighbourhood of
$\spec{H_c}$ and it is contained in the $\eps$-neigh\-bour\-hood
of $\nume{H_c}$. We will show later (corollary \ref{t16}) that if
$\Im(c)>0$, there exist an $E>0$ depending on $c$ such that
\[
   \pspe{\eps}{H_c} \subset \nume{H_c}
\]
for all $\eps<E$. Observe that if $c$ is real, since $H_c$ is
self-adjoint, $\pspe{\eps}{H_c}$ is actually equal to the
$\eps$-neighbourhood of $\spec{H_c}$, we stress that this property
is false in general for non-self-adjoint operators.

\begin{proposition}
\label{t19} For $\Im(c)>0$ fixed, the resolvent norm of $H_c$ is
symmetric with respect to the axis $c^{1/2}r,\, r\bel \R$. As a
consequence for all $\eps>0$ the \linebreak $\eps$-pseudospectrum
of $H_c$ is also symmetric with respect to this axis.
\end{proposition}

\Proof For all $a>0$ let
\[
   T_af(x):=a^{1/2}f(ax) \hspace{.3in} (x \bel \R),
\]
the operator $T_a$ is isometric in $\Lt(\R)$ and such that
$T_a^{-1}=T_{a^{-1}}$. Putting $\arg(c):=\tet$, since
\[
   T_{|c|^{1/4}} H_{\e^{i\tet}}  T_{|c|^{-1/4}} = |c|^{-1/2} H_c,
\]
it becomes evident that we can assume without lost of generality
$c=\e^{i\tet}$ for $0<\tet< \pi/2$.

For simplicity we rewrite the operator
\[
   H_c=P^2+cQ^2
\]
where $Pf(x):=if'(x)$ and $Qf(x)=xf(x)$ (respectively the quantum
mechanical observables of momentum and position). Thus, applying
Fourier transform we obtain
\begin{eqnarray*}
   \norm (P^2+cQ^2-z)^{-1} \norm & = & \norm (Q^2+cP^2-z)^{-1} \norm \\
   & = & \norm (P^2+\e^{-i\tet}Q^2-\e^{-i\tet}z)^{-1} \norm \\
   & = & \norm [(P^2+\e^{-i\tet} Q^2 - \e^ {-i\tet} z)^{-1}]^\ast \norm \\
   & = & \norm (P^2+cQ^2-\e^{i\tet}\overline{z})^{-1} \norm
\end{eqnarray*}
for all $z \not\in \spec{H_c}$, which is precisely our claim. \fin

\section{High Energy Eigenvalues}

In this section we show that if the coupling constant $c$ is such
that $\Im(c)>0$, and $z_\eta\bel \nume{H_c}$ parameterized by
$\eta>0$ is such that
\[
   |z_\eta - (b \eta +c \eta^p)|\to 0 \hspace{.3in} (\mathrm{as}\; \eta \to
   \infty)
\]
for some $b>0$ and $1/3 <p < 3$, then
\[
   \lim _{\eta \to \infty} \norm (H_c-z_\eta)^{-1} \norm = \infty.
\]
Such a result will be a consequence of theorem \ref{t18}. Notice
that by proposition \ref{t19} it is enough to assume that $1/3 < p
\leq 1$.

Our first aim is to obtain test functions $f_\eta \bel \Ic(\R)$,
parameterized by $\eta>0$, such that if
\[
   z_{\eta}= ic\eta^{1/2-\gam /2} +\alpha ^2 \eta ^{\gam} + c \alpha ^2 \eta ,
\]
where $\alpha >0$ and $1\leq \gam < 3$ are fixed constants
independent of $\eta>0$,
\begin{equation}
   \label{e10}
   \lim_{\eta \to \infty} \frac{\norm H_c f_{\eta} - z_{\eta} f_{\eta}
   \norm }{\norm f_{\eta} \norm } =0.
\end{equation}

We follow a similar procedure as the one given originally in
\cite{D2}. Let
\[
    \Phi (x) := \e ^{-\Psi (x)}
    \hspace{.5in} (x\bel \R),
\]
where the polynomial
\[
    \Psi (x_0 +s) := \Psi _1 s + \Psi _2 s^2 /2 + \Psi _3 s^3 /3
    \hspace{.5in} (s \bel \R),
\]
centered in $x_0:=\alpha \eta ^{1/2}$, has coefficients
\begin{eqnarray*}
     \Psi_1 & := & i \alpha \eta^{\gam /2} , \\
     \Psi_2 & := & -i c \eta ^{1/2 - \gam /2} , \\
     \Psi_3 & := & -\frac{ic}{2 \alpha} \eta^{-\gam /2} (1+c\eta^{1-\gam}) .\\
\end{eqnarray*}
Clearly
\[
    H_c \Phi (x_0 +s)  =  (p(s) + z_\eta ) \Phi (x_0 +s)
\]
where
\begin{equation}
    \label{e11}
    p(s)=c_1 s+c_2 s^2+c_3 s^3 +c_4 s^4
    \hspace{.5in} (s\bel \R)
\end{equation}
has coefficients
\begin{eqnarray*}
    c_1 & = & -\frac{i c}{\alpha} \eta ^{-\gam /2}(1 + c \eta ^{1-\gam} ), \\
    c_2 & = & 0, \\
    c_3 & = & \frac{c^2}{\alpha} \eta ^{1/2 - \gam} (1+c \eta ^{1-\gam}), \\
    c_4 & = & \frac{c^2}{4 \alpha ^2} \eta ^{-\gam}(1+ 2c\eta ^{1-\gam}
              +c^2 \eta ^{2-2 \gam}).
\end{eqnarray*}

Let us establish some properties of the function $\Phi$ above
defined. A straightforward calculation implies that
\[
    \left| \Phi (x_0 +s) \right| ^2 =  \exp \left[ -\beta _2 (\eta) s^2 -
                        \beta _3 (\eta) s^3 \right],
\]
for
\begin{eqnarray*}
    \beta _2 (\eta ) & := &  \Im (c) \eta ^{1/2 - \gam /2} >0 , \\
    \beta _3 (\eta ) & := &  \frac{\Im(c) (1+ 2 \Re (c) \eta ^{1-\gam})}
                             {3 \alpha} \eta^{-\gam /2} > 0.
\end{eqnarray*}
Taking derivatives of $| \Phi (x_0 +s) |^2$ with respect to $s$,
allows us to conclude that this function has a local maximum at
$s=0$ , and a local minimum at
\[
   s=s_0:= - \frac{2 \beta _2 (\eta)}{3 \beta _3 (\eta)}
          = -\frac{2 \alpha \eta ^{1/2}}{1+ 2 \Re (c) \eta ^{1-\gam }}
         \to -\infty
\]
as $\eta \to \infty$.

The required $f_\eta$ can be defined truncating the function
$\Phi$ as follows. It is routine to define a $\Ci (\R)$
compact support function
\[
    g(x)= \left\{
    \begin{array}{lll} 1 & {\rm if} & |x-x_0 | < \eta ^{\del _0} \\
    0 & {\rm if} & |x-x_0 | > 2 \eta ^{\del_0} \\
    \end{array} \right.
\]
for $\del _0 := \gam /6$, such that there exist constants $q_1
>0$ and $q_2 >0$ independent of $x,\, \eta ,\, \alpha$ or $\gam$,
with the property
\begin{eqnarray*}
    |g'(x)| & \leq & q_1 \eta ^{-\del_0} ,\\
    |g''(x)| & \leq & q_2 \eta ^{-2 \del_0} , \\
\end{eqnarray*}
for any $x\bel \R$. Then,
\[
    f_\eta (x):= g(x) \Phi (x)
    \hspace{.5in} (x \bel \R).
\]

The next two lemmas point out some properties of $f_\eta$ which
can be employed to demonstrate (\ref{e10}). The constants $a_k$
for $k=1,\,2,\, \ldots$ below are real and independent of $\eta$,
but can possibly depend on $\alpha$ or $\gam$.

\begin{lemma}
    \label{t14}
    For $f_\eta$ as above and $1\leq \gam <3$,
    there exist positive constants $a_1$, $a_2$ and $E_\gam$, independent of
    $\eta$, such that
    \[
        a_1 \eta ^{(\gam-1)/4} \leq \norm f_\eta \norm ^2
        \leq a_2 \eta ^{(\gam-1)/4}
    \]
    for all $\eta>E_\gam$.
\end{lemma}
\Proof Since $\del_0 <1/2$,
\[
    \lim _{\eta \to \infty} \frac{2 \beta _3 (\eta) \eta ^{3\del_0} t^3}
    {\beta_2 (\eta) \eta ^{2\del_0} t^2} = \lim_{\eta \to \infty}
    \frac{2 (1+2 \Re (c) \eta ^{1-\gam}) \eta ^{\del_0} t}
    {3 \alpha \eta ^{1/2}}
    = 0,
\]
uniformly for $0\leq t \leq 2$. Therefore, there exist $E_\gam>0$
such that for any $\eta> E_\gam$
\begin{equation}
   \label{e13}
   0 \leq \beta _3 (\eta) \eta ^{3 \del_0} t^3 \leq \frac{\beta _2 (\eta)
   \eta ^{2 \del_0}t^2}{2},
\end{equation}
for any $0 \leq t \leq 2$.

For the first inequality: If $\eta>E_\gam$, we have
\begin{eqnarray*}
    \norm f_\eta \norm ^2 & = & \int _{\R} |g(x_0 +s)|^2
               |\Phi (x_0 +s)|^2 ds \\
    & \geq & \int_0 ^{\eta ^{\del_0}} \exp [-\beta _2 (\eta) s^2
             -\beta _3 (\eta) s^3] ds \\
    & \geq &  \int_0 ^1 \exp [-\beta _2 (\eta) \eta ^{2\del_0} t^2
             -\beta _3 (\eta) \eta^{3\del_0} t^3] \eta ^{\del_0} dt. \\
    & \geq &  \int_0 ^1 \exp [-2 \beta _2 (\eta)
             \eta ^{2\del_0} t^2] \eta ^{\del_0} dt \\
    & \geq & \int_0 ^{\eta^{\del_0 +(1-\gam)/4}} \exp [-2 \Im (c) u^2]
             \eta ^{-(1-\gam)/4} du \\
    & \geq & a_1 \eta ^{-(1-\gam)/4}.
\end{eqnarray*}
Observe that we are using the fact that $\del_0 +
(1-\gam)/4=1/4-\gam /12>0$. The second inequality is similar. \fin

\gap

We now estimate the numerator in the left hand side of
(\ref{e10}). By (\ref{e11}) we have
\begin{eqnarray}
   \norm H_c f_\eta - z _\eta f_\eta \norm & \leq &
   \norm 2 g'(x_0 +s ) \Phi '(x_0 +s)
   \norm  + \norm g'' (x_0 +s ) \Phi (x_0 +s ) \norm  \nonumber\\
   & & + \sum_{k=1}^4 \norm c_k s^k  \Phi (x_0 +s )
   g(x_0 +s ) \norm \label{e12}.
\end{eqnarray}

\begin{lemma}
   \label{t17}
   Let $E_\gam>0$ as in lemma \ref{t14}. There exist
   positive constants $a_4$, $a_5$ and $a_6$, independent of $\eta>0$, such that
   \begin{eqnarray*}
      \norm 2 g'(x_0 +s ) \Phi '(x_0 +s ) \norm ^2 & \leq
      & a_4 \eta^{5\gam /6} \exp [-\frac{\Im(c)}{2}
      \eta^{(3-\gam)/6}] \\
      \norm g''(x_0 +s ) \Phi (x_0 +s ) \norm ^2 & \leq &
      a_5 \eta ^{-\gam/2} \exp[- \frac{\Im (c)}{2}
     \eta ^{(3-\gam)/6}] \\
     \norm s^k \Phi(x_0 +s) g(x_0 +s) \norm ^2 & \leq &
     a_6 \eta^{(\gam -1)(2k+1)/4}
   \end{eqnarray*}
   for all $\eta>E_\gam$ and $k=1,\,3$ or $4$.
\end{lemma}
\Proof We use similar arguments as the one provided in the proof
of lemma \ref{t14}.

Let
\[
   \Ome:=\{s \bel \R \, : \, \eta ^{\del_0} \leq |s| \leq 2 \eta ^{\del_0} \}.
\]
For all $s \bel \Ome$, we have
\begin{eqnarray*}
    |\Psi '(x_0 +s) |^2 & = & |\Psi _1 +\Psi _2 s +\Psi _3 s^2 |^2 \\
         & \leq & \eta ^{\gam}(\alp + 2 |c| \eta ^{\del_0 +1/2 - \gam}
                 +\frac{2 |c|}{\alp} \eta ^{2\del_0 - \gam} (1 + |c| \eta
                 ^{1-\gam}))^2 \\
         & \leq & \eta ^\gam a_3.
\end{eqnarray*}
Then for all $\eta>E_\gam$, by the conditions imposed on $g(x)$
above and by (\ref{e13}), we have
\begin{eqnarray*}
    \norm 2 g'(x_0 +s ) \Phi '(x_0 +s ) \norm ^2
        & = & 4 \int _\R |g'(x_0 +s)|^2 |\Psi '(x_0 +s)| ^2
          |\Phi (x_0 +s) |^2 ds \\
        & \leq & \frac{a_4}{2} \eta^{\gam -2\del_0} \int _\Ome \exp [-\beta _2 (\eta) s^2 -
          \beta _3 (\eta) s^3] ds \\
        & \leq & a_4 \eta^{\gam -2\del_0} \int _{\eta^{\del_0}}
                 ^{2\eta^{\del_0}} \exp [-\beta _2 (\eta) s^2 +
                 \beta _3 (\eta) s^3] ds \\
        & \leq & a_4 \eta^{\gam -\del_0} \int _{1}^{2}
                  \exp [-\frac{\beta _2 (\eta)}{2} \eta^{2\del_0} t^2] dt \\
        & \leq & a_4 \eta^{5\gam /6} \exp [-\frac{\Im(c)}{2}
                  \eta^{(3-\gam)/6}].
\end{eqnarray*}
The second estimate is similar so that we can find $a_5$ without
difficulty.

Finally for $\eta>E_\gam$ and $k=1,\,3$ or $4$, using again
(\ref{e13}),
\begin{eqnarray*}
   \norm s^k \Phi(x_0 +s) g(x_0 +s) \norm ^2
       & \leq &  \int _{-2 \eta ^{\del_0}}^{2 \eta ^{\del_0}} s^{2k}
                |\Phi (x_0 +s) |^2 ds \\
       & \leq & 2 \int _{0}^{2\eta^{\del_0}} s^{2k} \exp [-\beta _2 (\eta) s^2 +
                 \beta _3 (\eta) s^3] ds \\
       & \leq & 2 \eta^{(2k+1)\del_0} \int _{0}^{2} t^{2k} \exp [-\frac{\beta _2
                   (\eta)}{2} \eta^{2\del_0} t^2] dt \\
       & \leq & 2 \eta^{(\gam -1)(2k+1)/4} \int _{0}^{\infty}
                   u^{2k} \exp [-\frac{\Im (c)}{2} u^2] du \\
        & \leq & a_6 \eta^{(\gam -1)(2k+1)/4}. \find
\end{eqnarray*}

Notice that (\ref{e10}) can be easily obtained from lemma
\ref{t14}, lemma \ref{t17} and equation (\ref{e12}). Using such
result, we can achieve the following theorem.

\begin{theorem}
   \label{t18}
   Let $H_c$ be the complex harmonic oscillator such that $\Re(c)>0$
   and $\Im(c)>0$. If
   \[
      z_\eta= b \eta + c  \eta^{p}
   \]
   where $b > 0$ and $1/3 < p < 3$ are constants independent
   of $\eta>0$,
   then
   \[
      \lim_{\eta \to \infty} \norm (H_c - z_{\eta} )^{-1} \norm =
      \infty.
   \]
\end{theorem}
\Proof Assume $1/3<p\leq 1$. If we put the unitary operator $T_a$
for $a>0$ as in the proof of proposition \ref{t19}, recall that
for $\arg(c):=\tet$
\[
   T_{|c|^{1/4}} H_{\e^{i\tet}}  T_{|c|^{-1/4}} = |c|^{-1/2} H_c.
\]

For all $r>0$ and $\bet>0$ let
\[
   w_{\eta,r}:=\bet \eta + \bet r \e^{i\tet} \eta^p.
\]
Since
\begin{eqnarray*}
   \norm (H_{\e^{i\tet}}-w_{\eta,|c|})^{-1} \norm & = &
   \norm (|c|^{-1/2}H_c-w_{\eta,|c|})^{-1} \norm \\
   & = & |c|^{1/2} \norm (H_c-|c|^{1/2}w_{\eta,|c|})^{-1} \norm,
\end{eqnarray*}
using (1) putting $\alp^2=|c|^{1/2}\bet$ and thinking of $|c|=r$,
it is clear that
\[
   \lim_{\eta \to \infty} \norm (H_{\e^{i\tet}}-w_{\eta,r})^{-1}
   \norm = \infty
\]
for all $r>0$ and $\bet>0$.

Now if we put $\bet=|c|^{1/2}b$ and $r=|c|/b$, we obtain
\begin{eqnarray*}
   \norm (H_c - z_\eta)^{-1} \norm & = &
   \norm (|c|^{1/2} H_{e^{i\tet}} - z_\eta)^{-1} \norm \\
   & = & |c|^{-1/2} \norm (H_{e^{i\tet}} - w_{\eta,|c|/b})^{-1} \norm \to
   \infty
\end{eqnarray*}
as $\eta \to \infty$. The proof can be now completed by
proposition \ref{t19}. \fin

This theorem has some consequences for the $\eps$-pseudospectra of
$H_c$. It implies that for $\eps>0$ fixed, any curve $z_\eta\bel
\nume{H_c}$ parameterized by $\eta>0$ such that
\[
   |z_\eta - (b \eta +c \eta^p)|\to 0 \hspace{.3in} (\mathrm{as}\; \eta \to
   \infty)
\]
for some $b>0$ and $1/3 <p < 3$, will eventually be inside
$\pspe{\eps}{H_c}$ and it will stay there as $\eta \to \infty$. In
particular this shows that high energy eigenvalues are
increasingly unstable under small perturbations.

\section{Non-Self-Adjoint Mehler's Formula}

In this section we show that the bounded holomorphic semigroup of
contractions generated by $-H_c$ is compact in a maximal angular
sector. In order to deduce this property for the semigroup, we
obtain an explicit formula for the heat kernel of $H_c$ as in
\cite{D4} and show directly that this kernel is of Hilbert-Schmidt
type.

We say that for $\alp,\,\bet>0$, a parameterized family of bounded
operators $T_\tau$ in the Banach space $\BB$, is a {\it bounded
holomorphic semigroup of contractions} in the sector
$S(-\alp,\bet)$, if and only if:
\begin{enumerate}
   \item $T_{\tau_1+\tau_2}=T_{\tau_1}T_{\tau_2}$ for all $\tau_k\bel
   S(-\alp,\bet)$.
   \item $\norm T_\tau \norm \leq 1$ for all $\tau \bel
   S(-\alp,\bet)$.
   \item $T_\tau$ is a holomorphic family of operators in
   $\tau\bel S(-\alp,\bet)$.
   \item For all $f\bel \BB$ and $\eps >0$,
   \[
      \lim_{\tau \to 0} T_\tau f=f
   \]
   for $\tau$ inside $S(-\alp+\eps,\bet-\eps)$.
\end{enumerate}
It follows directly from the definition that, for
$-\alp<\tet<\bet$ fixed, $T_{\e^{i\tet}t}$ for $t>0$ is a $C_0$
one-parameter semigroup in the standard sense, \cite{ops}. The
{\it generator} of $T_\tau$ is, by definition, the infinitesimal
generator of the one-parameter semigroup $T_t$ where $t>0$.

For convenience, we shall put
\[
   S_c := \left\{ \begin{array}{lcr}
      \overline{S(-\frac{\pi}{2},\frac{\pi}{2}-\arg(c))}
      \setminus \{0\} & \mathrm{if} & \Im(c) > 0   \\
      S(-\frac{\pi}{2},\frac{\pi}{2})
      & \mathrm{if} & \Im(c) = 0   \\
      \overline{S(-\frac{\pi}{2}-\arg(c),\frac{\pi}{2})}
      \setminus \{0\} & \mathrm{if} & \Im(c) < 0.   \\
      \end{array} \right.
\]

\begin{theorem}
   \label{t2}
   Let $H_c$ the complex harmonic oscillator for $\Re(c)>0$ as defined above, then:
   \begin{enumerate}
   \item For $c$ fixed, $-H_c$ is the generator of a bounded holomorphic semigroup
   of contractions $\e^{-H_c\tau}$, with parameter $\tau$ in the open
   sector $\inte(S_c)$.
   \item If $\Im(c) > 0$, then  $iH_c$ and $-\e^{i(\pi/2-\arg(c))}H_c$
   are also generators of a one-parameter semigroup such that $\e^{-H_c\tau}$ is
   strong continuous for all $\tau \bel S_c$.
   \item If $\Im(c) < 0$, then $-iH_c$ and $-\e^{-i(\pi/2+\arg(c))}H_c$
   are also generators of a one-parameter
   semigroup such that $\e^{-H_c\tau}$ is strongly continuous for all $\tau \bel S_c$.
   \item For $\tau >0$ fixed, $\e^{-H_c\tau}$ is a holomorphic family of bounded
   operators parameterized by $c$ for all $\Re(c)>0$.
   \end{enumerate}
\end{theorem}

\Proof Properties 1 is deduced without difficulty from
\cite[theorem 1.24, p.492]{kat} and \cite[theorem 2.24]{ops}.

Property 2 and 3: if $\Im(c)>0$, operator $iH_c$ is maximal
dissipative operators, therefore is the generator of a
one-parameter semigroup (the same argument works for the other
three cases). From analyticity, follows strong continuity for all
$\tau \bel \inte(S_c)$ and strong continuity in the edges can be
checked using \cite[corollary 3.18]{ops}.

Observe that if $\Im(c)=0$, $\spec{\pm iH_c}$ is purely imaginary,
so we cannot apply the above argument.

Property 4: by the way in which we define its domain, the operator
$H_c$ is a holomorphic family of type (B) for $\Re(c)>0$. Using
the fact that holomorphic families of this type are locally
m-sectorial, this property can be demonstrated by analogy to
\cite[theorem 2.6, p.500]{kat}. \fin

\gap

For $\tau \bel S_c$ let the coefficients:
\begin{eqnarray*}
   \lam & := & \exp[-2c^{1/2}\tau], \\
   w_1 & := & c^{1/4} \lam^{1/2}[\pi (1-\lam^{2})]^{-1/2}, \\
   w_2 & := & \frac{c^{1/2}(1+\lam^2)}{2(1-\lam^2)}, \\
   w_3 & := & \frac{2c^{1/2}\lam}{(1-\lam^2)}.
\end{eqnarray*}
For all $x,\,y\, \bel \R$, put
\[
   K_c(\tau,x,y):= w_1 \exp \left[ w_3xy-w_2(x^2+y^2) \right],
\]
 and define the integral operator
\begin{equation}
   \label{e5}
   A_{c,\tau}f(x):=\int _\R K_c(\tau,x,y)f(y)dy,
\end{equation}
for all $f\bel \Lt (\R)$ for which this formula makes sense.

Fixing $\tau>0$, $K_c(\tau,x,y)$ is holomorphic for $\Re(c)>0$,
and fixing $c$ it is continuous for $\tau\bel S_c$ and analytic in
its interior. Observe that for $c$ real, $K_c(\tau,x,y)$ has poles
if $\tau$ runs along the complex axis; this situation is avoided
by the definition we have made of the angular region $S_c$.

As pointed out in \cite{D4}, modifications of the classical
argument (see \cite{cyc,ops}) show that
\begin{equation}
\label{e6}
   e^{-H_c\tau}=A_{c,\tau}
\end{equation}
for all $c>0$ and $\tau>0$. Using analytic continuation this
equality can be extended to complex $c$ and $\tau$.

\begin{theorem}[Non-self-adjoint Mehler's Formula]
   \label{t3}
   Let $K_c(\tau,x,y)$ be defined as above.
   Then for all $\Re(c)>0$ and $\tau \bel S_c$,
   \[
      e^{-H_c\tau}f(x)=\int _\R K_c(\tau,x,y)f(y)dy \hspace{.3in}
      (f\bel \Lt(\R)).
   \]
\end{theorem}

The rest of this section is devoted to proving this theorem. First
we establish some local bounds for $|K_c(\tau,x,y)|$, this allows
us to show that $A_{c,\tau}$ is a holomorphic family of
Hilbert-Schmidt operators in both parameters and then we obtain
the desired equality by analytic continuation.

\begin{lemma}
   \label{t6}
   For all $\tau_0 \bel S_c$, there exists a neighbourhood
   $\tau_0\bel V_0 \subset \C$ and real constants $\alp_1$, $\alp_2$ and $\alp_3$
   such that for all $\tau\bel V_0 \cap S_c$:
   $$
      |K_c(\tau,x,y)| \leq \alp_1 \exp\left[\alp_3xy-\alp_2(x^2+y^2)\right] \hspace{.5in}
      (\mathrm{all}\, x,\,y \bel \R).
   $$
   The constants satisfy: $\alp_1,\,\alp_2>0$ and $2\alp_2 \pm \alp_3>0$,
   locally uniformly in $c$ and $\tau$.
\end{lemma}
\Proof Clearly
\[
   |K_c(\tau,x,y)| = |w_1| \exp[\Re(w_3)xy-\Re(w_2)(x^2+y^2)].
\]
Since $w_i$ are continuous in $\tau$ and $c$, it is enough to show
that $\Re(w_2)>0$ and $\Re(2w_2 \pm w_3)>0$, for all $\Re(c)>0$
and $\tau\bel S_c$.

Let $\tet_c:=\arg(c^{1/2})$ and $\tet_\tau:=\arg(\tau)$. Without
loss of generality we will assume that $\Im(c)>0$ and
$|c^{1/2}|=2$. Define the Moebius transforms by
\[
   M_{\pm}(z):=\frac{1\pm z}{1 \mp z} \hspace{.5in} (z\bel \C).
\]
Then
\[
   w_2=\frac{c^{1/2}}{2} M_+ (\lam^2)=\e^{i\tet_c} M_+ (\lam^2)
\]
and
\[
   2w_2 \pm w_3 = c^{1/2}\frac{1\pm \lam}{1 \mp \lam}=2 \e^{i\tet_c} M_\pm (\lam).
\]

Notice that $M_+$ maps the interior of the unitary disk into
the open right half plane, and multiplying by $\e^{i\tet_c}$
rotates about the origin by $\tet_c$. Every disk centered in
the origin, of radius $r_1<1$, changes under $M_+$ into a disk
$D_2$ of radius
\[
   r_2:=\frac{2r_1}{1-r_1^2}
\]
centered at
\[
   c_2:=\frac{1+r_1^2}{1-r_1^2} >0.
\]
If $r_1$ is close to 0, $r_2$ is small, $c_2$ is close to 1 and
the rotation of $D_2$ about the origin by $\tet_c$ remains in the
right half plane. If $r_1$ is close to 1, $r_2$ is big, and rotating by
$\tet_c$ can send some points of $D_2$ to the left half plane. In spite
of this
possibility, notice that points in $D_2$ which are to the right of
$c_2$ do not cross to the left half plane when rotated. It is elementary to
find the maximum radius $r_1$ which allows $\e^{i\tet_c}D_2$ stay
in the right half plane. We call, the disk with center the origin having this
radius, the {\it critical disk}.

Now if we consider $\tet_c$ and $\tet_\tau$ fixed,
$\lam^2$ describes a spiral with radius decreasing
exponentially as $|\tau|$ increases, starting in 1 and coiling
about the origin. We leave to the reader, to check that all these
spirals cross the critical disk with real part large enough to
guarantee that after the mapping by $M_+$ and the rotation by
$\tet_c$, the resulting curve stays in the right half plane. This gives
$\Re(w_2)>0$ and a similar analysis gives $\Re(2w_2 \pm w_3)>0$.
\fin

\pagebreak

\begin{lemma}
    \label{t7}
    For all $\Re(c)>0$ and $\tau\bel S_c$, the operator $A_{c,\tau}$ as defined in
    (\ref{e5})
    is of Hilbert-Schmidt type and
    \begin{enumerate}
        \item For $c$ fixed, $A_{c,\tau}$ is norm continuous
        for $\tau \bel S_c$ and holomorphic for $\tau \bel \inte(S_c)$.
        \item For $\tau>0$ fixed, $A_{c,\tau}$ is holomorphic family of
        bounded operators for $\Re(c)>0$.
    \end{enumerate}
\end{lemma}

\Proof Lemma \ref{t6} implies that for all $\tau_0 \bel S_{c}$,
there exist a neighbourhood $\tau_0 \bel V_0 \subset \C$ and a
constant $M<\infty$ such that:
\[
    \int \int |K_c(\tau,x,y)|^2 dx dy < M,
\]
for all $\tau\bel V_0\cap S_c$, where the constant $M$ depends
locally uniformly on $c$. Therefore $A_{c,\tau}$ is
Hilbert-Schmidt operator and its norm depends locally uniformly on
$c$ and $\tau$.

By the dominated convergence theorem, it is elementary to check
that
\[
    \norm K_c(\tau,\cdot,\cdot)-K_c(\tau_0,\cdot,\cdot) \norm_2 \to
    0
\]
as $\tau\to \tau_0$. This provides norm continuity of $A_{c,\tau}$ in
$\tau\bel S_c$ for $c$ fixed. Finally the analyticity of the
family of operator in both variables, can be deduced without
difficulty by differentiating under the integral
sign. \fin

\gap

By theorem \ref{t2} and lemma \ref{t7}, we know that the families
of bounded operators $\e^{-H_c\tau}$ and $A_{c,\tau}$ are both
holomorphic in each parameter. They coincide when those parameters
are real. If we fix $\tau>0$, by analytic continuation in $c$,
equation (\ref{e6}) is also true for all $\Re(c)>0$. Now fixing
$c$ with $\Re(c)>0$, by analytic continuation in $\tau$, the same
equation is true for all $\tau \bel \inte(S_c)$. Finally, since
both families are strongly continuous in $\tau$ at the edges of
$S_c$, we have proved theorem \ref{t3} as requested.

Lemma \ref{t7} is one of the crucial points for the analysis of
pseudospectra that we intend to carry out in the next section. In
particular, we use strongly the fact that $\e^{-H_c\tau}$ is a
compact operator for all $\tau \bel S_c$.

\section{Pseudospectra}
Our aim in this section is to employ the compactness of the
bounded holomorphic semigroup generated by $-H_c$ (see previous
section) to obtain estimates on the resolvent norm of $H_c$ inside
$S(0,\arg(c))$. The technique is based on the spectral radius
formula for the semigroup and a convenient Jordan decomposition of
the operators involved. For simplicity we will assume from now on
that $\Im(c)>0$, but with a few corrections, the same results and
techniques apply to the case $\Im(c)<0$.

We establish first a formula for the spectral radius of
one-parameter semigroups; a proof of this fact can be found in
\cite[theorem 1.22]{ops}. If $\e^{-Tt}$ is a $C_0$ one-parameter
semigroup, the limit
\[
   a:=\lim_{t\to \infty} t^{-1} \log \norm \e^{-Tt} \norm
\]
always exists with $-\infty \leq a < \infty$, and the spectral
radius
\[
   \srad{\e^{-Tt}}:=\max \left\{|\lam| \,:\, \lam \bel \spec{\e^{-Tt}}
   \right\} =\e^{at}
\]
for all $t>0$. This implies in particular that for all
$\alp<-a$
\begin{equation}
\label{e7}
   \lim_{t\to \infty} e^{\alp t} \norm \e^{-Tt} \norm  = 0.
\end{equation}

Applying this estimate to $\e^{-H_c\tau}$ for $\tau$ on the edges
of $S_c$, we can show that the resolvent norm of $H_c$ is
uniformly bounded in lines parallel and close enough to the edges
of $S(0,\arg(c))$.

\begin{theorem}
   \label{t11}
   For fixed $\Im(c)>0$, consider the complex  parameters $z_{\rm low}=\eta + i \eps$
   and
   $z_{\rm upp}=c(\eta-i \eps)/|c|$, with $\eta>0$ and $\eps>0$.
   Then, for $0<d<\Im(\lam_0)$ there exists a constant $M_{c,d}<\infty$,
   independent of
   $\eta$ and
   $\eps$, such
   that
   \[
      \norm (H_c-z_{\rm low})^{-1} \norm \leq M_{c,d}
   \]
   and
   \[
      \norm (H_c-z_{\rm upp})^{-1} \norm \leq M_{c,d},
   \]
   for all $\eta>0$ and $0\leq \eps \leq d$.
\end{theorem}

\Proof By theorem \ref{t3} and lemma \ref{t7}, the bounded
holomorphic semigroup $\e^{-H_c\tau}$ is compact for all $\tau$ in
the maximal sector $S_c$. Therefore, \cite{ops},
\[
   \spec{\e^{-H_c\tau}}=\{0 \} \cup \{\e^{-\lam_n \tau } \, :\,
   n=0,\,1,\, \ldots \}.
\]
By proposition \ref{t19} we can just concentrate on the lower
edge. Putting $\tau=-it$ for $t>0$, we obtain
\[
   \srad{\e^{iH_ct}}=\e^{-\Im(\lam_0)t}.
\]
Fix $0<d<\Im(\lam_0)$; by formula (\ref{e7}) with $T=-iH_c$,
$a=-\Im(\lam_0)$ and \linebreak $\alp=(d-a)/2$, there exists
$t_\alp>0$ such that
\[
   \norm \e^{iH_ct} \norm \leq \e^{-\alp t}
\]
for all $t>t_\alp$. Then
\begin{eqnarray*}
   \norm (H_c - (\eta+i\eps) )^{-1} \norm & = &
       \norm (i\eta-\eps -iH_c )^{-1} \norm \\
   & \leq & \int _0^\infty \e^{\eps s} \norm \e^{iH_c s} \norm ds  \\
   & \leq & \int _0^{\tau_{\alp}} \e^{\eps s} ds +
      \int _{\tau_{\alp}}^{\infty} \e^{\left(\eps-\alp\right) s} ds \\
   & \leq &  \int _0^{\tau_{\alp}} \e^{\eps s} ds +
      \int _{\tau_{\alp}}^{\infty} \e^{\left(d-\alp \right) s} ds  \\
   & \leq &  \int _0^{\tau_{\alp}} \e^{\eps s} ds +
      \int _{0}^{\infty} \e^{-(\Im(\lam_0)-d)s/2}
       ds  \\
   & \leq & M_{c,d}. \find
\end{eqnarray*}

This result provides information about the shape of the
$\eps$-pseudospectra of operator $H_c$. The result below will be
extended in theorem \ref{t15}.

\begin{corollary}
\label{t16}
   For all $0<\del<1$ there exists an $\eps>0$ such that
   \[
      \pspe{\eps}{H_c} \subset S(0,\arg(c)) + \del c^{1/2}.
   \]
\end{corollary}

As a consequence of this corollary together with proposition
\ref{t20}, we obtain
\[
   \pspe{\eps}{H_c} \subset \nume{H_c}
\]
for $\eps$ small enough.

\gap

The situation for the other eigenvalues, requires more
care and the argument involves a Jordan decomposition of the problem.
Let $Q_n$ the spectral projector of operator $H_c$ associated to
the eigenvalue $\lam_n$, {\it i.e.}
\[
   Q_n f:=\frac{1}{2\pi i} \int_{\gam_n} (z-H_c)^{-1} f dz
\]
where $\gam_n$ is a smooth curve whose interior just contains
eigenvalue $\lam_n$. Observe that $Q_n$ are not orthogonal
projections in $\Lt(\R)$.

For $m=0,\,1,\, \ldots$ put
\[
   P_m:=\sum_{n=0}^m Q_n.
\]
It is clear that $P_m$ is the spectral projector with rank $m$
associated to eigenvalues $\lam_0,\,\ldots,\, \lam_m$. Notice that
$\Lt(\R)$ can be decomposed as the direct sum of the closed
subspaces $\rang{Q_n}$ for $n\leq m$ and $\rang{I-P_m}$, in the
sense that the subspaces are linearly independent and
\[
   \Lt(\R)= \rang{Q_0}+ \ldots +\rang{Q_m} +\rang{I-P_m}.
\]

It is easy to show that each of the $m+2$ subspaces above is
invariant under the semigroup $\e^{-H_c\tau}$ for all $\tau \bel
S_c$. Even more, the generator of the bounded holomorphic
semigroups of contractions
\[
   \e^{-H_c \tau}|_{\rang{Q_n}}
\]
and
\[
   \e^{-H_c \tau}|_{\rang{I-P_m}},
\]
are respectively $-H_c|_{\rang{Q_n}}$ and $-H_c|_{\rang{I-P_m}}$.
Notice as well that, by compactness of the restriction of the
semigroup, we have
\begin{equation}
    \label{e9}
    \spec{\e^{-H_c \tau}|_{\rang{I-P_m}}} = \{0\} \cup \{\e^{-\lam_n \tau}
    \}_{n=m+1}^{\infty}.
\end{equation}

In order to extend corollary \ref{t16} beyond the first eigenvalue,
observe that for all $z \not\in \spec{H_c}$
\[
   (H_c-z)^{-1} = \sum_{n=0}^m
   (H_c-z)^{-1}Q_n +
   (H_c-z)^{-1}(I-P_m).
\]
Then if we call $\kap_m:=1+ \sum_{n=0}^m \norm Q_n \norm$, a
straightforward calculation allows us to obtain the following
estimate: for all $z \not\in \spec{H_c}$
\begin{equation}
   \label{e8}
   \norm (H_c-z)^{-1} \norm \leq \kap_m \left( \sum_{n=0}^m \norm
   (H_c|_{\rang{Q_n}}-z)^{-1}\norm + \norm
   (H_c|_{\rang{I-P_m}}-z)^{-1}\norm \right).
\end{equation}
With the help of this inequality, we can achieve the theorem
below.

\begin{theorem}
   \label{t15}
   Let $H_c$ the complex harmonic oscillator such
   that $\Re(c)>0$ and $\Im(c)>0$. For all $0<\del<1$ and
   $m=0,1,\ldots$ there exists an $\eps>0$ such that
   \[
      \pspe{\eps}{H_c} \subset \left[ S(0,\arg(c)) +
      (\lam_{m+1}-\del) \right] \cup \bigcup_{n=0}^m
      \{z\bel\C \,:\, |z-\lam_n|<\del\}.
   \]
\end{theorem}
\Proof To estimate the first sum of norms in (\ref{e8}), since
$\rang{Q_n}$ is the one dimensional subspace generated by the
$n^{th}$ eigenvector $\Psi_n$, the operator $H_c|_{\rang{Q_n}}$
acts on this subspace as the operator of multiplication by
$\lam_n$. Then for all $z\neq \lam_n$,
\[
   \norm (H_c|_{\rang{Q_n}}-z)^{-1} \norm = \frac{1}{|\lam_n-z|}.
\]

To estimate the last resolvent norm in (\ref{e8}), we use equation
(\ref{e9}), applying an analogous of theorem \ref{t11} to
$H_c|_{\rang{I-P_m}}$ instead of $H_c$. Notice that now the first
eigenvalue is $\lam_{m+1}$. \fin

\gap

In the notation of Aslanyan and Davies, \cite{AD1,D4},
\[
   \kap_m = 1+\sum_{n=1}^{m} \kap(\lam_m),
\]
where $\kap(\lam_n)$ is the index of instability of the eigenvalue
$\lam_n$. Based on the results in \cite{D4}, as $n$ increases, the
indices of instability $\kap(\lam_n)$ grow faster than any power
of $n$. This is reflected in the above theorem in the fact that as
$m$ gets bigger the $\eps$ we must choose gets exponentially
smaller.

This theorem confirms the numerical calculations made by Davies in
\cite{D2}, where he uses a computer package to find level curves
for the resolvent norm of a discretization of the operator $H_c$.

Going beyond the scope of theorem \ref{t15}, our conjecture is as
follows. Let \linebreak $0< p <1/3$, $m=0,1,\ldots$ and $0<\del<1$
be fixed. Let $b_{m,p}
>0$ such that there exist $E>0$ (possibly depending on $m$ or $p$)
verifying
\[
   b_{k,p}E+cE^p=\lam_{m}.
\]
Put
\[
   z_\eta:=b_{m,p} \eta + c\eta^{p} \hspace{.3in} (\eta>0)
\]
and
\[
   \Ome_{m,p}:= \left\{ |z_\eta|\e^{i\tet} \in \C \,:\, \eta \geq E,\,
   \arg(z_\eta)\leq \tet \leq \arg(c\overline{z_{\eta}}/|c|)
   \right\}.
\]

\begin{conjecture}
There exists an $\eps>0$ such that
\[
   \pspe{\eps}{H_c} \subset
   \Ome_{m,p} \cup \bigcup_{n=0}^m
      \{z\bel\C \,:\, |z-\lam_n|<\del\}.
\]
\end{conjecture}

This would be a substantial improvement of the results provided in
this paper. The case $p=0$ is precisely theorem \ref{t15}. Because
of theorem \ref{t18} the statement is false for $1/3< p \leq 1$ so
in this sense the constraint $0<p<1/3$ is optimal.

\gap

{\bf Acknowledgments} The author would like to thank to Prof. E.
B. Davies for valuable discussions of different aspects of the
paper and to Dr. A. Aslanyan for helpful suggestions. This work
was supported by ``Fundaci\'on Gran Mariscal de Ayacucho'',
Venezuela, E-211-1357-1997-1.

\vspace{1in}

\setlength{\baselineskip}{13pt}

\setlength{\parindent}{0cm}
\vspace{.5in}

Department of Mathematics
\newline King's College
\newline Strand (The Strand)\newline London WC2R 2LS \newline England \newline

e-mail: lboulton@mth.kcl.ac.uk
\end{document}